# Semistable Abelian Varieties over $\mathbb{Z}[\frac{1}{6}]$ and $\mathbb{Z}[\frac{1}{10}]$

Frank Calegari

### Abstract

We prove there do not exist any non-zero semistable Abelian varieties over $\mathbb{Z}[\frac{1}{6}]$ or over $\mathbb{Z}[\frac{1}{10}]$. Our results are contingent on the GRH discriminant bounds of Odlyzko. Combined with recent results of Brumer–Kramer and of Schoof, this result is best possible: it shows that semistable Abelian varieties over $\mathbb{Z}[1/N]$, where $N$ is square-free, exist precisely for $N \notin \{1, 2, 3, 5, 6, 7, 10, 13\}$.

## 1   Introduction.

In 1985, Fontaine [3] proved a conjecture of Shafarevich to the effect that there do not exist any (non-zero) Abelian varieties over $\mathbb{Z}$ (equivalently, Abelian varieties $A/\mathbb{Q}$ with good reduction everywhere). Fontaine's approach was via finite group schemes over local fields. In particular, he proved the following theorem:

**Theorem 1.1 (Fontaine)** *Let $G_\ell$ be a finite flat group scheme over $\mathbb{Z}_\ell$ killed by $\ell$. Let $L = \mathbb{Q}_\ell(G_\ell) := \mathbb{Q}_\ell(G_\ell(\overline{\mathbb{Q}}_\ell))$. Then*

$$v(\mathfrak{D}_{L/\mathbb{Q}_\ell}) < 1 + \frac{1}{\ell - 1}$$

*where $v$ is the valuation on $L$ such that $v(\ell) = 1$, and $\mathfrak{D}_{L/\mathbb{Q}_\ell}$ is the different of $L/\mathbb{Q}_\ell$.*

If $G_\ell$ is the restriction of some finite flat group scheme $G/\mathbb{Z}$ killed by $\ell$ then $\mathbb{Q}(G)$ is *a fortiori* unramified at primes outside $\ell$. In this context, the result of Fontaine is striking since it implies that the field $\mathbb{Q}(G)$ has particularly small root discriminant. If $A/\mathbb{Q}$ has good reduction everywhere, then it has a smooth proper Néron model $\mathcal{A}/\mathbb{Z}$, and $G := \mathcal{A}[\ell]/\mathbb{Z}$ is a finite flat group scheme. Using the discriminant bounds of Odlyzko [7], Fontaine showed that for certain small primes $\ell$, for every $n$, either $A/\mathbb{Z}$ or some isogenous variety has a rational $\ell^n$-torsion point for every $n$. Reducing $A$ modulo $p$ for some prime $p$ of good reduction (in this case, any prime), one finds Abelian varieties (of fixed dimension $d$) over $\mathbb{F}_p$ with at least $\ell^n$ rational points. One knows, however, that isogenous Abelian varieties over $\mathbb{F}_p$ have an equal and thus bounded number of points. This contradiction proves Fontaine's Theorem.

If one considers Abelian varieties $A/\mathbb{Q}$ such that $A$ has good reduction outside a single prime $p$, one can no longer expect non-existence results. Indeed, there exist Abelian varieties with good reduction everywhere except at $p$. One such class of examples are the Jacobians of modular curves $X_0(p^n)$, which have positive genus for every $p$ and sufficiently large $n$. A natural subclass of Abelian varieties, however, are the semistable ones. By considering the modular Abelian varieties $J_0(N)$, with $N$ squarefree, one finds non-zero semistable Abelian varieties unramified outside $N$ for all $N \notin \{1, 2, 3, 5, 6, 7, 10, 13\}$. A reasonable conjecture to make is that there are no semistable Abelian varieties over $\mathbb{Z}[1/N]$ for $N$ in this set. Fontaine's Theorem is the case $N = 1$. Recently Brumer and Kramer [1] prove this result for $N \in \{2, 3, 5, 7\}$, and (by quite different methods) Schoof [9] for $N \in \{2, 3, 5, 7, 13\}$. In this paper, we treat the remaining cases $N \in \{6, 10\}$. Since we shall exploit results from both Brumer-Kramer [1] and Schoof [9], we briefly recall the main ideas now.

Schoof's approach is similar in spirit to Fontaine's. Instead of working with finite flat group schemes over $\mathbb{Z}$, one considers finite flat group schemes over $\mathbb{Z}[\frac{1}{p}]$, where $p$ is prime. In order to avoid group schemes arising from non-semistable Abelian varieties, one uses the following fact due to Grothendieck ([4], Exposé IX, Proposition 3.5):

**Theorem 1.2 (Grothendieck)** *Let $A$ be an Abelian variety with semistable reduction at $p$. Then the action of inertia at $p$ on the $\ell^n$-division points of $A$ is rank two unipotent; i.e., as an endomorphism, for $\sigma \in \mathcal{I}_p$,*

$$(\sigma - 1)^2 A[\ell^n] = 0.$$

*In particular, $\mathcal{I}_p$ acts through its maximal pro-$\ell$ quotient, which is procyclic.*

Thus one may restrict attention to finite flat group schemes $G/\mathbb{Z}[\frac{1}{p}]$ of $\ell$-power order such that inertia at $p$ acts through its maximal pro-$\ell$ quotient. The key step of Schoof's approach is to show that any such group scheme admits a filtration by the group schemes $\mathbb{Z}/\ell\mathbb{Z}$ and $\mu_\ell$. Using this filtration, along with various extension results (in the spirit of Mazur [6], in particular Proposition 2.1 pg. 49 and Proposition 4.1 pg. 58) for group schemes over $\mathbb{Z}[\frac{1}{p}]$, one shows as in Fontaine that for each $n$, some variety isogenous to $A$ has rational torsion points of order $\ell^n$.

The approach of Brumer and Kramer is quite different. Although, as in Schoof and Fontaine, they use discriminant bounds to control $\mathbb{Q}(A[\ell])$ for particular $\ell$, they seek a contradiction not to any local bounds but to a theorem of Faltings. Namely, they construct infinitely many pairwise non isomorphic but isogenous varieties, contracting the finiteness of this set (as follows from Faltings [2], Satz 6, pg. 363). The essential difference in the two approaches, however, is that Brumer and Kramer use the explicit description of the Tate module $\mathbb{T}_\ell$ of $A$ at a prime $p$ of semistable reduction. Such a description is once more due to Grothendieck [4]. Both of these approaches fail (at least naïvely) to work when $N = 6$ or $10$. Using Schoof's approach, one runs into a problem (when $N = 6$, for example) because $\mu_5$ admits many non-isomorphic finite flat group scheme extensions by $\mathbb{Z}/5\mathbb{Z}$ over $\mathbb{Z}[\frac{1}{6}]$, whereas no non-trivial extensions exist over either $\mathbb{Z}[\frac{1}{2}]$ or $\mathbb{Z}[\frac{1}{3}]$. Using Brumer and Kramer's approach, one difficulty that arises is that the field $\mathbb{Q}(A[5])$ fails to have a unique prime above the bad primes 2 or 3, as fortuitously happens in the cases they consider. We do, however, use a key theorem from Brumer and Kramer's paper, and so in the next section we recall some of their definitions and results.

## 1.1 Notation.

Let $p \in \mathbb{Z}$ be a prime number. Let $\mathcal{D}_p = \mathrm{Gal}(\overline{\mathbb{Q}}_p/\mathbb{Q}_p)$ denote the local Galois group at $p$. For a Galois extension of global fields $L/\mathbb{Q}$, we denote a decomposition group at $p$ by $\mathcal{D}_p(L/\mathbb{Q})$. This is well defined up to conjugation, or equivalently, up to an embedding $\overline{\mathbb{Q}} \hookrightarrow \overline{\mathbb{Q}}_p$ which we shall fix when necessary. In the same spirit, let $\mathcal{I}_p = \mathrm{Gal}(\overline{\mathbb{Q}}_p/\mathbb{Q}_p^{unr})$, and let $\mathcal{I}_p(L/\mathbb{Q})$ be an inertia group at $p$ as a subgroup of $\mathcal{D}_p(L/\mathbb{Q})$ and of $\mathrm{Gal}(L/\mathbb{Q})$. One notes that $\mathcal{I}_p$ is normal in $\mathcal{D}_p$. Let $\mathcal{M}$ be a $\mathcal{D}_p$ module, $\overline{\mathcal{M}}$ a $\mathcal{D}_p$ module killed by $\ell$ for some $\ell \neq p$, and $\widehat{\mathcal{M}}$ a $\mathrm{Gal}(\overline{\mathbb{Q}}/\mathbb{Q})$ module, also killed by $\ell$. A "finite" group scheme $G/R$ will always mean a group scheme $G$ finite and flat over Spec $R$.

# 2 Local Considerations.

## 2.1 Preliminaries.

In this section we introduce some notation and results from the paper of Brumer and Kramer [1].

Let $A/\mathbb{Q}$ be an Abelian variety of dimension $d > 0$ with semistable reduction at $p$. Let $\ell$ be a prime different from $p$, and consider the Tate module $\mathbb{T}_\ell(A/\mathbb{Q}_p)$. Let $\mathcal{M}_1(p) = \mathbb{T}_\ell(A/\mathbb{Q}_p)^{\mathcal{I}}$, and let $\mathcal{M}_2(p)$ be the subspace of $\mathbb{T}_\ell(A/\mathbb{Q}_p)$ orthogonal to $\mathcal{M}_1(p)(\hat{A})$ under the Weil paring:

$$e_\infty : \mathbb{T}_\ell(A) \times \mathbb{T}_\ell(\hat{A}) \longrightarrow \mathbb{Z}_\ell(1).$$

Since $A$ is semistable, there exist inclusions:

$$0 \subseteq \mathcal{M}_2(p) \subseteq \mathcal{M}_1(p) \subseteq \mathbb{T}_\ell(A/\mathbb{Q}_p).$$

Since $\mathcal{I}_p$ is normal in $\mathcal{D}_p$, $\mathcal{M}_1(p)$ and $\mathcal{M}_2(p)$ are $\mathrm{Gal}(\overline{\mathbb{Q}}_p/\mathbb{Q}_p)$ modules. Let $\mathcal{A}/\mathbb{Z}$ be a Néron model for $A$. Let $\mathcal{A}^0_{\mathbb{F}_p}$ be the connected component of the special fibre of $\mathcal{A}$ at $p$. It is an extension of an Abelian variety of dimension $a_p$ by a torus of dimension $t_p = d - a_p$. One has $\dim(\mathcal{M}_2(p)) = t_p$ and $\dim(\mathcal{M}_1(p)) = t_p + 2a_p = d + a_p$.

**Definition 2.1 (Brumer-Kramer)** *Let $i(A, \ell, p)$ denote the minimal integer $n \geq 1$ such that $\mathbb{Q}_p(A[\ell^n])$ is ramified at $p$. Call $i(A, \ell, p)$ the "effective stage of inertia". Since $A$ has bad reduction at $\ell$, $i(A, \ell, p)$ is finite by the criterion of Néron–Ogg–Shafarevich.*

Let $\Phi_A(p) = \mathcal{A}/\mathcal{A}^0(\overline{\mathbb{F}}_p)$ be the component group of $A$ at $p$. Recall the following result from [1]:

**Theorem 2.1 (Brumer-Kramer)** *Let $\overline{\mathcal{M}}_1(p)$ and $\overline{\mathcal{M}}_2(p)$ denote the projections of $\mathcal{M}_1(p)$ and $\mathcal{M}_2(p)$ to $A[\ell]$. Suppose that $\kappa$ is a proper $\mathrm{Gal}(\overline{\mathbb{Q}}/\mathbb{Q})$ submodule of $A[\ell]$ and let $\phi : A \longrightarrow A'$ be the $\mathbb{Q}$-isogeny with kernel $\kappa$. Then*

$$\mathrm{ord}_\ell(\Phi_{\hat{A}'}(p)) - \mathrm{ord}_\ell(\Phi_{\hat{A}}(p)) = \dim\ \kappa \cap \overline{\mathcal{M}}_2(p) + \dim\ \kappa \cap \overline{\mathcal{M}}_1(p) - \dim\ \kappa.$$

*Moreover, if $\overline{\mathcal{M}}_2(p) \subseteq \kappa \subseteq \overline{\mathcal{M}}_1(p)$, then $i(A', \ell, p) = i(A, \ell, p) + 1$.*

Brumer and Kramer use this theorem to construct infinitely many non-isomorphic varieties isogenous to $A$. This contradicts Faltings' Theorem. Although we shall also use Faltings' Theorem, our final contradiction will come from showing that $A$ (or some isogenous variety) has too many points over some finite field, contradicting Weil's Riemann hypothesis, much as in the approach of Schoof [9].

## 2.2 Results.

Our main results are the following:

**Theorem 2.2** *Let $A/\mathbb{Q}$ be an Abelian variety with semistable reduction, and good reduction outside $2$ and $3$. Assuming the GRH discriminant bounds of Odlyzko, $A$ has dimension $0$.*

**Theorem 2.3** *Let $A/\mathbb{Q}$ be an Abelian variety with semistable reduction, and good reduction outside $2$ and $5$. Assuming the GRH discriminant bounds of Odlyzko, $A$ has dimension $0$.*

The use of the GRH is impossible to avoid using our approach. The proof of Theorem 2.3 is very similar to the proof of Theorem 2.2, although some additional complications arise. Thus we restrict ourselves first to the case $N = 6$, and then later explain how our proof can be adapted to work for $N = 10$. One main ingredient is the following result, proved in section 3:

**Theorem 2.4** *Let $G/\mathbb{Z}[\frac{1}{6}]$ be a finite group scheme of $5$-power order such that inertia at $2$ and $3$ acts through a procyclic $5$-group. Then $G$ has a filtration by the group schemes $\mathbb{Z}/5\mathbb{Z}$ and $\mu_5$. Moreover, if $G$ is killed by $5$, then $\mathbb{Q}(G) \subseteq K$, where $K := \mathbb{Q}(\sqrt[5]{2}, \sqrt[5]{3}, \zeta_5)$.*

In particular, if $A/\mathbb{Q}$ is a semistable Abelian variety with good reduction outside $2$ and $3$, and $\mathcal{A}/\mathbb{Z}$ is its Néron model, then for each $n$ the finite group scheme $\mathcal{A}[5^n]/\mathbb{Z}[\frac{1}{6}]$ has a filtration by the group schemes $\mathbb{Z}/5\mathbb{Z}$ and $\mu_5$. Moreover, $\mathbb{Q}(A[5]) \subseteq K$. This result (and its proof) is of the same flavour as results in Schoof [9]. One such result from that paper we use is the following (a special case of Theorem 3.3 and the proof of corollary 3.4 in *loc. cit.*):

**Theorem 2.5 (Schoof)** *Let $p = 2$ or $3$. Let $G/\mathbb{Z}[\frac{1}{p}]$ be a finite group scheme of $5$-power order such that inertia at $p$ acts through a procyclic $5$-group. Then $G$ has a filtration by the group schemes $\mathbb{Z}/5\mathbb{Z}$ and $\mu_5$. Moreover, the extension group* $\mathrm{Ext}^1(\mu_5, \mathbb{Z}/5\mathbb{Z})$ of group schemes over $\mathbb{Z}[\frac{1}{p}]$ is trivial, and there exists an exact sequence of group schemes:

$$0 \longrightarrow M \longrightarrow G \longrightarrow C \longrightarrow 0$$

where $M$ is a diagonalizable group scheme over $\mathbb{Z}[\frac{1}{p}]$, and $C$ is a constant group scheme.

In sections 2.3, 2.4 and 2.5 we shall assume there exists a semistable Abelian variety $A/\mathbb{Z}[\frac{1}{6}]$, and derive a contradiction using Theorem 2.4.

## 2.3 Construction of Galois Submodules.

The proof of Brumer and Kramer relies on the fact that for Abelian varieties with semistable reduction at one prime $p \in \{2, 3, 5, 7\}$, there exists an $\ell$ such that there is a unique prime above $p$ in $\mathbb{Q}(A[\ell])$. In this case, the $\mathcal{D}_p$ modules $\overline{\mathcal{M}}_1(p)$ and $\overline{\mathcal{M}}_2(p)$ are automatically $\mathrm{Gal}(\overline{\mathbb{Q}}/\mathbb{Q})$ modules, and so one has a source of $\mathrm{Gal}(\overline{\mathbb{Q}}/\mathbb{Q})$ modules with which to apply Theorem 2.1. This approach fails in our case (at least if $\ell = 5$) since Theorem 2.4 allows the possibility that $\mathbb{Q}(A[5])$ could be as big as $K := \mathbb{Q}(2^{1/5}, 3^{1/5}, \zeta_5)$, and 2 and 3 split into 5 distinct primes in $\mathcal{O}_K$. On the other hand, something fortuitous does happen, and that is that the inertia subgroups $\mathcal{I}_p(K/\mathbb{Q})$ for $p = 2$, 3 are *normal* subgroups of $\mathrm{Gal}(K/\mathbb{Q})$, when *a priori* they are only normal subgroups of $\mathcal{D}_p(L/\mathbb{Q})$. Using this fact we may construct global Galois modules from the local $\mathcal{D}_p$ modules $\overline{\mathcal{M}}_1(p)$ as follows.

**Lemma 2.1** *Let $F = \mathbb{Q}(A[\ell])$, $G = \mathrm{Gal}(F/\mathbb{Q})$, and $H \subseteq G$ be a normal subgroup of $G$. Let $\overline{\mathcal{M}}$ be a subgroup of $A[\ell]$ fixed pointwise by $H$. Let $\widehat{\mathcal{M}}$ be the $\mathrm{Gal}(\overline{\mathbb{Q}}/\mathbb{Q})$ submodule generated by the points of $\overline{\mathcal{M}}$. Then $\mathbb{Q}(\widehat{\mathcal{M}}) \subseteq E$, where $E$ is the fixed field of $H$.*

**Proof.** By Galois, it suffices to show that $\widehat{\mathcal{M}}$ is fixed by $H$. Any sum or multiple of elements fixed by $H$ is clearly fixed by $H$. Thus it remains to show that any Galois conjugate $P^g$ with $g \in G$ and $P \in \overline{\mathcal{M}}$ is also fixed by $H$. For this we observe that

$$(P^g)^h = (P^{ghg^{-1}})^g = P^g$$

Since $ghg^{-1} \in H$. □

Throughout, let $\widehat{\mathcal{M}}_1(p)$ be the $\mathrm{Gal}(\overline{\mathbb{Q}}/\mathbb{Q})$ module generated $\overline{\mathcal{M}}_1(p)$, considered as a subgroup of $A[\ell]$ after choosing some embedding $\overline{\mathbb{Q}} \hookrightarrow \overline{\mathbb{Q}}_p$ (this definition depend upon the embedding, but this ambiguity does not cause any problems). From Lemma 2.1, $\widehat{\mathcal{M}}_1(p)$ is fixed by $\mathcal{I}_p(K/\mathbb{Q})$ and so

$$\mathbb{Q}(\widehat{\mathcal{M}}_1(2)) \subseteq \mathbb{Q}(\zeta_5, 3^{1/5}), \qquad \mathbb{Q}(\widehat{\mathcal{M}}_1(3)) \subseteq \mathbb{Q}(\zeta_5, 2^{1/5}).$$

We now apply Theorem 2.1 with $\kappa = \widehat{\mathcal{M}}_1(2)$. Let $A' = A/\kappa$. Then

$$\mathrm{ord}_5(\Phi_{\hat{A}'}(2)) - \mathrm{ord}_5(\Phi_{\hat{A}}(2)) = \dim\ \kappa \cap \overline{\mathcal{M}}_2(2) + \dim\ \kappa \cap \overline{\mathcal{M}}_1(2) - \dim\ \kappa.$$

Since by construction $\overline{\mathcal{M}}_2(2) \subseteq \overline{\mathcal{M}}_1(2) \subseteq \kappa$, this quantity equals $2d - \dim \kappa \geq 0$. In particular, $A$ can not be isomorphic to $A'$ unless $\kappa = A[5]$. Thus by Faltings' Theorem, after a finite number of isogenies $A'[5] = \kappa = \widehat{\mathcal{M}}_1(2)$ and $A'[5]$ is unramified at 2. Replace $A$ by $A'$. Since $A[5]$ is unramified at 2, $\mathcal{A}[5]$ prolongs to a finite group scheme over $\mathbb{Z}[\frac{1}{3}]$. From Theorem 2.5, we infer that there exists an exact sequence of group schemes:

$$0 \longrightarrow \mu_5^m \longrightarrow A[5] \longrightarrow (\mathbb{Z}/5\mathbb{Z})^n \longrightarrow 0$$

where $m + n = 2d$.

**Lemma 2.2** *In the sequence above, $m = n = d$. $A$ has ordinary reduction at 5.*

**Proof**. The Néron model of $A' = A/\mu_5^m$ contains the group scheme $(\mathbb{Z}/5\mathbb{Z})^n$. Specializing to the fibre over $\mathbb{F}_5$ we find that:

$$(\mathbb{Z}/5\mathbb{Z})^n \hookrightarrow A'_{\mathbb{F}_5}[5].$$

The $p$-rank of the $p$-torsion subgroup of an Abelian variety in characteristic $p$ is at most the dimension $d$, with equality only if $A$ is ordinary at $p$. Thus $n \leq d$. Applying the same argument to $\hat{A}$ we find that $m \leq d$ and thus $n = m = d$, and $A$ has ordinary reduction at 5. $\square$

Thus we may assume for any $A$ with $\mathrm{ord}_5(\Phi_{\hat{A}}(2))$ maximal (or, by a similar argument $\mathrm{ord}_5(\Phi_{\hat{A}}(3))$ maximal) there exists an exact sequence of $\mathrm{Gal}(\overline{\mathbb{Q}}/\mathbb{Q})$ modules:

$$0 \longrightarrow \mu_5^d \longrightarrow A[5] \longrightarrow (\mathbb{Z}/5\mathbb{Z})^d \longrightarrow 0$$

We now divide our proof by contradiction into two cases. In the first case we assume that $A$ has mixed reduction at 2 or at 3. In the second case we assume that $A$ has purely toric reduction at both 2 and 3.

## 2.4 $A$ has Mixed Reduction at $2$ or $3$.

Let $\mathrm{ord}_5(\Phi_{\hat{A}}(2))$ be maximal. Then from Lemma 2.2 there is an exact sequence:

$$0 \longrightarrow \mu_5^d \longrightarrow A[5] \longrightarrow (\mathbb{Z}/5\mathbb{Z})^d \longrightarrow 0.$$

If $A$ has mixed reduction at 2 then $a_2 > 0$, and $\mathcal{M}_1(2)$ has dimension $t_2 + 2a_2 = d + a_2 > d$. In particular, $\kappa := \overline{\mathcal{M}}_1(2) \cap \mu_5^d$ is non-trivial and defines a diagonalizable $\mathrm{Gal}(\overline{\mathbb{Q}}/\mathbb{Q})$ submodule of $A[5]$. We now apply Theorem 2.1. Let $A' = A/\kappa$. We find that

$$\mathrm{ord}_5(\Phi_{\hat{A}'}(2)) - \mathrm{ord}_5(\Phi_{\hat{A}}(2)) = \dim \kappa \cap \overline{\mathcal{M}}_2(2) + \dim \kappa \cap \overline{\mathcal{M}}_1(2) - \dim \kappa.$$

Since $\kappa \subseteq \overline{\mathcal{M}}_1(2)$, the last two terms cancel, and $\mathrm{ord}_5(\Phi_{\hat{A}'}(2))$ is also maximal. Hence we may repeat this process, thereby constructing morphisms $A \longrightarrow A^{(n)}$ with larger and larger kernels $\kappa_n$, where $\kappa_n$ has a filtration by $\mu_5$'s.

**Lemma 2.3** *Any extension of diagonalizable group schemes of* 5*-power order over* $\mathbb{Z}[\frac{1}{6}]$ *is diagonalizable.*

**Proof.** By taking Cartier duals, it suffices to prove the dual statement for constant group schemes: Any extension of 5-power order constant group schemes over $\mathbb{Z}[\frac{1}{6}]$ is constant. Any extension of $\mathbb{Z}/5\mathbb{Z}$ by $\mathbb{Z}/5\mathbb{Z}$ over $\mathbb{Z}[\frac{1}{6}]$ is defined over a 5-extension of $\mathbb{Q}$, unramified over 6. From class field theory, since $\mathbb{Z}$ is an integral domain, such extensions are classified by $(\mathbb{Z}/6\mathbb{Z})^*$. Since this group has order coprime to 5, this proves the claim. $\square$

For all $n$, there exist exact sequences

$$0 \longrightarrow \kappa_n \longrightarrow A[5^{k(n)}] \longrightarrow M \longrightarrow 0.$$

The variety $\hat{A}/M^\vee$ contains the arbitrarily large constant group scheme $\kappa_n^\vee$. This contradicts the uniform boundedness of the number of points locally for all varieties isogenous to $\hat{A}$.

If $A$ does not have purely toric reduction at 3, a similar argument applies.

## 2.5 $A$ has Purely Toric Reduction at $2$ and $3$.

Under this assumption, for $p = 2$ or 3, $\mathcal{M}_2(p) = \mathcal{M}_1(p)$, and so we write both as $\mathcal{M}(p)$. Again we assume that $\text{ord}_5(\Phi_{\hat{A}}(2))$ is maximal. In particular, we may assume that $\widehat{\mathcal{M}}(2) = A[5]$, that $A[5]$ is defined over $\mathbb{Q}(\zeta_5, 3^{1/5})$, and that we have an exact sequence

$$0 \longrightarrow \mu_5^d \longrightarrow A[5] \longrightarrow (\mathbb{Z}/5\mathbb{Z})^d \longrightarrow 0.$$

**Lemma 2.4** $\widehat{\mathcal{M}}(3) = \mu_5^d$.

**Proof.** Fix an embedding $\overline{\mathbb{Q}} \hookrightarrow \overline{\mathbb{Q}}_2$ such that the image of $3^{1/5}$ lands in $\mathbb{Q}_2$. First we show that $\overline{\mathcal{M}}(2) \cap \mu_5^d = \{0\}$. If not, then $\overline{\mathcal{M}}(2)$ would not surject onto $(\mathbb{Z}/5\mathbb{Z})^d$, and the elements of $\overline{\mathcal{M}}(2)$ could not possibly generate[1] $A[5]$ as a $\text{Gal}(K/\mathbb{Q})$ module. Thus by dimension considerations, as a $\mathbb{F}_5$ vector space, $A[5] = \mu_5^d \oplus \overline{\mathcal{M}}(2)$.

$\overline{\mathcal{M}}(2)$ is a $\mathcal{D}_2$ module. Consider generators $\{P_1, \ldots P_d\}$ for $\overline{\mathcal{M}}(2)$. Let $G = \text{Gal}(\mathbb{Q}(\zeta_5, 3^{1/5})/\mathbb{Q})$. Then for our embedding of $\overline{\mathbb{Q}}$ into $\overline{\mathbb{Q}}_2$,

$$D_2 := \mathcal{D}_2((\mathbb{Q}(\zeta_5, 3^{1/5})/\mathbb{Q}) \simeq \text{Gal}(\mathbb{Q}(\zeta_5, 3^{1/5})/\mathbb{Q}(3^{1/5})) \simeq \mathbb{Z}/4\mathbb{Z} = \{\tau\}.$$

Since $\overline{\mathcal{M}}(2)$ is a $D_2$ module, the $P_i$ are permuted by elements of $D_2$. Thus we may write the Galois action of $G = \{\sigma, \tau | \sigma^5 = \tau^4 = 1, \tau\sigma\tau^{-1} = \sigma^2\}$ on our basis $A[5] = \mu_5^d \oplus \overline{\mathcal{M}}(2)$ as follows, where $M_d$ and $N_d$ are $d \times d$ matrices, $\text{Id}_d$ is the identity matrix, and $\chi$ is the cyclotomic character.

$$\tau := \begin{pmatrix} \chi \cdot \text{Id}_d & 0 \\ 0 & M_d \end{pmatrix}, \qquad \sigma := \begin{pmatrix} \text{Id}_d & N_d \\ 0 & \text{Id}_d \end{pmatrix}.$$

---
[1]Another way to eliminate this possibility is as follows: if there was some intersection, we could apply the same argument as section 2.4 where such an intersection was useful.

Since any subgroup of $\mu_5^d$ is a $\text{Gal}(\overline{\mathbb{Q}}/\mathbb{Q})$ submodule, $\overline{\mathcal{M}}(2)$ can only generate $A[5]$ if $N_d$ is surjective. Thus $N_d$ is invertible.

Let us now consider the situation locally at 3. The decomposition group at 3 is the entire Galois group $G$, and the inertia group $\mathcal{I}_3$ is equal to $\{\sigma\} \subset G$. We show that $\mu_5^d \subset \widehat{\mathcal{M}}(3)$ and $\widehat{\mathcal{M}}(3) \subset \mu_5^d$. Since $\mathcal{I}_3$ acts faithfully on $\overline{\mathcal{M}}(2) \subset A[5]$, and since $\widehat{\mathcal{M}}(3)$ is unramified at 3, $\widehat{\mathcal{M}}(3) \cap \overline{\mathcal{M}}(2) = \{0\}$ and so $\widehat{\mathcal{M}}(3) \subseteq \mu_5^d$. On the other hand, $|\widehat{\mathcal{M}}(3)| \geq |\overline{\mathcal{M}}(3)| = 5^d = |\mu_5^d|$. Thus we are done. $\square$

We now apply Theorem 2.1 again with $\kappa = \widehat{\mathcal{M}}(3) = \mu_5^d$. If $A' = A/\mu_5^d$, then since $\overline{\mathcal{M}}(3) = \widehat{\mathcal{M}}(3)$, $i(A', 5, 3) = i(A, 5, 3) + 1 \geq 2$. On the other hand, we see from the exact sequence for $A[5]$ that $(\mathbb{Z}/5\mathbb{Z})^d \subset A'[5]$. From Theorem 2.4 and Lemma 2.2 we infer that there exists an exact sequence:

$$0 \longrightarrow (\mathbb{Z}/5\mathbb{Z})^d \longrightarrow A'[5] \longrightarrow \mu_5^d \longrightarrow 0.$$

Replace $A$ by $A'$.

**Lemma 2.5** *$A[5]$ is defined over $\mathbb{Q}(\zeta_5)$. There is only one prime above 3 in the extension $\mathbb{Q}(A[5])/\mathbb{Q}$.*

**Proof**. Consider the action of $\mathcal{I}_5$ on $A[5]$. Since $A$ is ordinary at 5, $A[5]$ (as an $\mathcal{I}_5$ module) is the extension of a constant module of rank $d$ by a cyclotomic module of rank $d$. The $(\mathbb{Z}/5\mathbb{Z})^d$ inside $A[5]$ must intersect trivially with the cyclotomic module. Thus it provides a splitting of $A[5]$ as an $\mathcal{I}_5$ module into a product of a cyclotomic module and a constant module. Thus $\mathbb{Q}(A[5])$ is unramified over $\mathbb{Q}(\zeta_5)$. The maximal extension of $\mathbb{Q}(\zeta_5)$ inside $K$ unramified at $1 - \zeta_5$ is $\mathbb{Q}(\zeta_5, 18^{1/5})$. Since $\mathbb{Q}(A[5])$ is also unramified over 3 (as $i(A, 5, 3) \geq 2$), $\mathbb{Q}(A[5])$ must be exactly $\mathbb{Q}(\zeta_5)$. The second statement is clear. $\square$

Thus, as in [1], $\kappa = \widehat{\mathcal{M}}(3) = \overline{\mathcal{M}}(3)$ is a $\text{Gal}(\overline{\mathbb{Q}}/\mathbb{Q})$ module, and applying Theorem 2.1 once more with $A' = A/\kappa$, since $\kappa = \overline{\mathcal{M}}(3)$ we find that

$$i(A', 5, 3) = i(A, 5, 3) + 1 \geq 3.$$

Replace $A$ by $A'$. In particular, $A[5^2]$ is unramified at 3. Thus by Theorem 2.5 there exists a filtration:
$$0 \longrightarrow M \longrightarrow A[5^2] \longrightarrow C \longrightarrow 0$$
where $M$ is a diagonalizable group scheme, and $C$ is a constant group scheme. Let $q \in \mathbb{Z}$ be a prime of good reduction. We observe that the varieties $A/M$ and $\hat{A}/C^\vee$ contain constant subgroup schemes of order $\#C$ and $\#M$ respectively. It follows from Weil's Riemann Hypothesis that Abelian varieties of dimension $d$ over $\mathbb{F}_q$ have at most $(1 + \sqrt{q})^{2g}$ points. Thus

$$5^{4g} = \#A[5^2] = \#C \#M \leq (1 + \sqrt{q})^{4g}.$$

Choosing $q = 7$, say, then since $5 > 1 + \sqrt{7}$, we have a contradiction. This completes the proof of Theorem 2.2 up to Theorem 2.4, which we prove now.

# 3   Group Schemes over $\mathbb{Z}[1/6]$.

First, some preliminary remarks on group schemes. Here we follow Schoof [9].

Let $(\ell, N) = 1$. Let $\underline{C}$ be the category of finite group schemes $G$ over $\mathbb{Z}[1/N]$ satisfying the following properties:

1. $G = G[\ell]$.
2. The action of $\sigma \in \mathcal{I}_p$ on $G(\overline{\mathbb{Q}}_p)$ is either trivial or cyclic of order $\ell$.

For example, $\mathbb{Z}/\ell\mathbb{Z}$ and $\mu_\ell$ are objects of $\underline{C}$. As remarked in [9], this category is closed under direct products, flat subgroups and flat quotients. Thus, to prove that any object of $\underline{C}$ has a filtration by $\mathbb{Z}/\ell\mathbb{Z}$ and $\mu_\ell$ it suffices to show that the only simple objects of $\underline{C}$ are $\mathbb{Z}/\ell\mathbb{Z}$ and $\mu_\ell$. If $A/\mathbb{Z}[1/N]$ is a semistable Abelian variety, then from Theorem 1.2, $\mathcal{A}[\ell] \in \underline{C}$. Another class of examples are due to Katz–Mazur ([5] Chapter 8, Interlude 8.7, [9]). For any unit $\epsilon \in \mathbb{Z}[1/N]$ there is a corresponding group scheme $G_\epsilon$ of order $\ell^2$ killed by $\ell$. It is an extension of $\mathbb{Z}/\ell\mathbb{Z}$ by $\mu_\ell$, and is defined over $\mathbb{Q}(\zeta_\ell, \epsilon^{1/\ell})$.

Let $N = 6$ and $\ell = 5$. To prove that the only simple objects of $\underline{C}$ are $\mu_5$ and $\mathbb{Z}/5\mathbb{Z}$, it suffices to show that any object of $\underline{C}$ is defined over the field $K$, where $K = \mathbb{Q}(\zeta_5, 2^{1/5}, 3^{1/5})$, because of the following result:

**Lemma 3.1** *Let $G/\mathbb{Z}[1/N]$ be a simple group scheme killed by $\ell$, where $(N, \ell) = 1$. Let $L = \mathbb{Q}(G(\overline{\mathbb{Q}}))$ and suppose that $\mathrm{Gal}(L/\mathbb{Q}(\zeta_\ell))$ is an $\ell$-group. Then $G$ is either $\mathbb{Z}/\ell\mathbb{Z}$ or $\mu_\ell$.*

**Proof.** Since any $\ell$-group acting on $(\mathbb{Z}/\ell\mathbb{Z})^d$ has at least one (in fact $\ell - 1$) non-trivial fixed points, there exists a point $P$ of $G$ defined over $\mathbb{Q}(\zeta_\ell)$. Since $G$ is simple, $P$ generates $G$ as a Galois module and thus $\mathbb{Q}(G) \subseteq \mathbb{Q}(\zeta_\ell)$. Since $(N, \ell) = 1$, and since $G$ is unramified outside $\ell$, $G$ prolongs to a finite group scheme over $\mathbb{Z}$, killed by $\ell$, and defined over $\mathbb{Q}(\zeta_\ell)$. Since the $(\ell - 1)^{th}$ roots of unity are in $\mathbb{F}_\ell^*$, any simple subgroup scheme of $G$ has order $\ell$. From Oort–Tate [8], the finite group schemes of order $p$ over $\mathbb{Z}$ are $\mathbb{Z}/p\mathbb{Z}$ and $\mu_p$. $\square$

Let $G$ be an object of $\underline{C}$. To prove that $\mathbb{Q}(G) \subseteq K$ it clearly suffices to prove the same inclusion for any group scheme which contains $G$ as a direct factor. Consider the field $L = \mathbb{Q}(G \times G_{-1} \times G_2 \times G_3)$. One sees (from the definition of $G_\epsilon$) that $K := \mathbb{Q}(\zeta_5, 2^{1/5}, 3^{1/5}) \subseteq L$. We prove that $L = K$. Using the estimates of Fontaine [3] we obtain an upper bound on the ramification of $L$ at 5. Since inertia at 2 and 3 acts through a cyclic subgroup of order 5, we also have ramification bounds at 2 and 3. As in Schoof [9] and Brumer–Kramer [1], we obtain has the following estimate of the root discriminant:

$$\delta_L < 5^{1 + \frac{1}{5-1}} 2^{1 - \frac{1}{5}} 3^{1 - \frac{1}{5}} = 5^{5/4} 6^{4/5} = 31.349 < 31.645.$$

From the discriminant bounds of Odlyzko [7], under the assumption of GRH, one concludes that $[L : \mathbb{Q}] < 2400$ and thus $[L : K] < 24$. In particular, $L/\mathbb{Q}$ is a solvable extension, and thus we can apply tools from class field theory.

**Remark**. Without the GRH, we are unable to bound $[L : \mathbb{Q}]$ since 31 exceeds the limits of current unconditional discriminant bounds.

Our calculations in this section could be shortened by more reliance on computer calculation. However, for exposition we include as much class field theory as we can do by hand. This leads us to consider several group theory lemmas which allow us to do computations in smaller fields.

The root discriminant of $K$ is $\delta_K = 5^{23/20}6^{4/5}$, and so $L/K$ is unramified outside primes above 5. Let $F = \mathbb{Q}(\zeta_5, 576^{1/5}) = \mathbb{Q}(\zeta_5, 18^{1/5}) = \mathbb{Q}(\zeta_5, 24^{1/5})$. Then $F/\mathbb{Q}(\zeta_5)$ is unramified at $1 - \zeta_5$. The prime $1 - \zeta_5$ splits completely in $F$, as

$$5 = \pi_1 \ldots \pi_5, \qquad \pi_i = \left(5, \left(\frac{\sqrt[5]{576} - 1}{1 - \zeta_5}\right) - i\right) \qquad N_{F/\mathbb{Q}}(\pi_i) = 5.$$

The extension $K/F$ is totally ramified at all primes $\pi_i$, $\pi_i = \mathfrak{p}_i^5$ for all $i$, and $N_{K/\mathbb{Q}}(\mathfrak{p}_i) = 5$. Let us consider the factorization of $\mathfrak{p}_i$ in $L$. Since $L/\mathbb{Q}$ is Galois, the ramification exponents are equal for all $i$. Thus we may write

$$\mathfrak{p}_i = \prod_{j=1}^{r_{L/K}} \mathfrak{P}_{i,j}^{e_{L/K}}, \qquad N_{L/\mathbb{Q}}(\mathfrak{P}_{i,j}) = 5^{f_{L/K}}, \qquad r_{L/K}e_{L/K}f_{L/K} = [L:K].$$

## 3.1  $L/K$ Tame.

In this section we assume that $L/K$ is a *tame* extension, of order coprime to 5, and prove that $L = K$.

**Lemma 3.2** $[L : K] < 10$.

**Proof**. Since $L/K$ is tame, $\mathfrak{D}_{L/K} = \prod_{j=1}^{r_{L/K}} \mathfrak{P}_{i,j}^{e_{L/K}-1}$, where $\mathfrak{D}_{L/K}$ is the different. Thus:

$$\Delta_{L/K} = N_{L/K}(\mathfrak{D}_{L/K}) = \prod_{i=1}^{5} \mathfrak{p}^{r_{L/K}f_{L/K}(e_{L/K}-1)}.$$

Since $N_{K/\mathbb{Q}}(\mathfrak{p}_i) = 5$, $\operatorname{ord}_5(N_{L/K}(\Delta_{L/K})) = 5[L:K](1 - 1/e_{L/K}) < 5[L:K]$. Using the transitivity property of the discriminant [10] we find:

$$\delta_L = \delta_K \cdot N_{K/\mathbb{Q}}(\Delta_{L/K})^{1/[L:\mathbb{Q}]} < \delta_K \cdot 5^{5/[K:\mathbb{Q}]} = 5^{23/20}6^{2/3}5^{5/100} = 28.925.$$

This contradicts the Odlyzko bounds [7]. If $[L : \mathbb{Q}] \geq 1000$, then assuming the GRH, $\delta_L > 29.094$. Thus $[L : K] < 10$. $\square$

**Remark**. If $5 \mid [L : K]$ and $[L : K] < 10$ then $[L : K] = 5$ and $\operatorname{Gal}(L/K)$ is tame if and only if it is unramified. We shall consider the case $L/K$ unramified in section 3.3. We may therefore assume that $[L : K]$ has order coprime to 5.

**Lemma 3.3** *Let $G$ be a finite group, and let $G' = [G, G]$ be its commutator subgroup. Suppose moreover that $G/G' \cong \mathbb{Z}/5\mathbb{Z}$, and that $|G'| < 10$. Then $G' = \{1\}$.*

**Proof**. It suffices to note that for all groups $G'$ of order less than 10, $|\text{Aut}(G')|$ is coprime to 5. □

**Lemma 3.4** *If $L/K$ is a tame extension of degree coprime to 5, then $L = K$.*

**Proof**. Let $H$ be the field $\mathbb{Q}(\zeta_5, 2^{1/5})$. We have the following exact sequence of groups:

$$0 \longrightarrow \text{Gal}(L/K) \longrightarrow \text{Gal}(L/H) \longrightarrow \mathbb{Z}/5\mathbb{Z} \longrightarrow 0$$

By Lemma 3.3, either $L = K$, or $\text{Gal}(L/K)$ is not the commutator subgroup of $\text{Gal}(L/H)$. Thus since $[L:K]$ has order coprime to 5, $\text{Gal}(L/H)^{ab}$ is not a 5-group. Hence $H$ admits a Galois extension $E/H$ contained in $L$, not contained in $K$, and of order coprime to 5.

**Sub-lemma 1** *$E/H$ is unramified at 2 and 3.*

**Proof**. Let $p \in \{2,3\}$. Consider ramification degrees $e_p$. One has

$$e_p(E/H) \mid e_p(L/H) = e_p(L/K)e_p(K/H).$$

Moreover, $L/K$ is unramified at primes above 2 and 3, and $[K:H] = 5$. Thus $e_p(E/\mathbb{Q}(\zeta_5))$ is a power of 5, which must be 1, since $5 \nmid [E:H]$. □

**(Continuation of Lemma)** Thus $E/H$ is a non-trivial Abelian extension of degree coprime to 5 and unramified outside $\pi_H$. Such extensions are classified by class field theory. One has by `pari` that $\text{Cl}(\mathcal{O}_H) = 1$. On the other hand, $H/\mathbb{Q}$ is totally ramified at 5, and so $(\mathcal{O}_H/\pi_5 \mathcal{O}_H)^* \simeq \mathbb{F}_5^*$ which is generated by the global units $(1+\sqrt{5})/2 \equiv -2$ and $-1$. Thus $E$ does not exist. This proves that $L = K$. □

## 3.2 $L/K$ Wild.

In this section, we assume that $L/K$ is wildly ramified and of degree 10, 15 or 20.

**Lemma 3.5** *Let $H$ be a group of order 10, 15 or 20. Let $G$ be an extension of $\mathbb{Z}/5\mathbb{Z}$ by $H$. Then $G^{ab}$ is not a 5-group.*

**Proof**. Let $H'$ be the 5-Sylow subgroup of $H$. Then since $5(1+5) > 20$, $H'$ is normal. Thus we have the following exact sequence:

$$0 \longrightarrow H' \longrightarrow H \longrightarrow H'' \longrightarrow 0.$$

Since $H''$ is Abelian the commutator subgroup of $H$ is a subgroup of $H''$. To show that $G^{ab}$ is not a 5-group it suffices to show that the commutator subgroup of $G$ also lies within the 5-Sylow subgroup of $H$. Let $\tau$ be an element of $G$ that maps to a generator of $\mathbb{Z}/5\mathbb{Z}$. The action of conjugation by $\tau$ on $H$ is via an automorphism of degree 5. To show that $[\tau, h] \in H'$ it suffices to show that for any automorphism $\sigma$ of degree 5 on $H$, $\sigma(h)h^{-1} \in H'$. Since all elements of order 5 lie in $H'$, $H'$ is preserved by $\sigma$. Yet $\text{Aut}(\mathbb{Z}/5\mathbb{Z}) \simeq \mathbb{Z}/4\mathbb{Z}$, and thus $\sigma$ fixes $H'$. Thus $\sigma$ maps to an element of $\text{Aut}(H'')$. Since $\text{Aut}(H'')$ has order coprime to 5 for $|H''| \leq 4$, $\sigma$ also acts trivially on the quotient. Hence for any automorphism $\sigma$ of degree 5 on $H$, $\sigma(h)h^{-1} \in H'$ and we are done. $\square$

Let $H$ be the field $\mathbb{Q}(\zeta_5, 2^{1/5})$. We have the following exact sequence of groups:

$$0 \longrightarrow \mathrm{Gal}(L/K) \longrightarrow \mathrm{Gal}(L/H) \longrightarrow \mathrm{Gal}(K/H) \longrightarrow 0.$$

By Lemma 3.5, $\mathrm{Gal}(L/H)^{ab}$ is not a 5-group. Thus $H$ admits an Abelian extension of degree coprime to 5. The non-existence of such an extension was proved in Lemma 3.4.

## 3.3 $L/K$ of degree 5.

Finally, it remains to show that $L/K$ is not wildly ramified of degree 5, or unramified over $K$. Assume otherwise. $\mathrm{Gal}(L/\mathbb{Q}(\zeta_5))$ is a group of order 125 that surjects onto $\mathbb{Z}/5\mathbb{Z} \oplus \mathbb{Z}/5\mathbb{Z}$. There are three groups up to isomorphism with this property. All of them admit at least one morphism to $\mathbb{Z}/5\mathbb{Z}$ with kernel $\mathbb{Z}/5\mathbb{Z} \oplus \mathbb{Z}/5\mathbb{Z}$ that factor through the map to $\mathrm{Gal}(K/\mathbb{Q}(\zeta_5))$. Thus there exists a field $E/\mathbb{Q}(\zeta_5)$, contained within $K$, such that $\mathrm{Gal}(L/E) \simeq \mathbb{Z}/5\mathbb{Z} \oplus \mathbb{Z}/5\mathbb{Z}$.

**Lemma 3.6** *There exists an intermediate field $L/F/E$ such that $F$ is not equal to $K$ and $F/E$ is unramified at primes above 2 and 3.*

**Proof**. Since the root discriminant of $L$ locally at 2 and 3 is bounded by $2^{4/5}$ and $3^{4/5}$ respectively, this lemma is obvious if the root discriminant for $E$ obtains these bounds, since then any subgroup of $\mathrm{Gal}(L/E) = \mathbb{Z}/5\mathbb{Z} \oplus \mathbb{Z}/5\mathbb{Z}$ not corresponding to $K$ will produce the required $F$. Thus we may assume that $E = \mathbb{Q}(p^{1/5}, \zeta_5)$ with $p$ equal to 2 or 3. Assume at $p = 2$. Since $K/E$ is ramified at primes above 3, it suffices to find an $F \subset L$ unramified at primes above 3. The tame ramification group $\mathcal{I}_3(L/E)$ is of order 5, since by considering 3 exponents of the root discriminant,

$$\delta_{L,3} = N_{E/\mathbb{Q}}(\Delta_{L/E})^{1/[L:\mathbb{Q}]} = 3^{1-1/e_{L/K}} \leq 3^{4/5}.$$

Thus we see that the fixed field $F$ of $\mathcal{I}_3(L/E) \subset \mathrm{Gal}(L/E)$ *is* unramified at 3 above $E$. Moreover, $F$ is not $K$ since $K/E$ is ramified at 3. An identical argument works for $p = 3$. $\square$

**Lemma 3.7** *If $E/\mathbb{Q}$ is wildly ramified at 5 then either $F/E$ is unramified at 5 or $\Delta_{F/E} = \pi_E^8$ where $\pi_E$ is the unique prime above 5 in $E$. If $E = \mathbb{Q}(\zeta_5, 24^{1/5})$ then $\Delta_{F/E}$ divides $(\pi_{E,1} \ldots \pi_{E,5})^8$, where $\pi_{E,i}$ are the primes above 5.*

**Proof**. Suppose that $E/\mathbb{Q}$ is wildly ramified. We may assume that $F/E$ is also wildly ramified, since otherwise it is unramified, and we are done. Suppose that $N_{E/\mathbb{Q}}(\Delta_{F/E}) \geq 5^{10}$. Then

$$\delta_{F,5} = \delta_{E,5} N_{E/\mathbb{Q}}(\Delta_{F/E})^{1/100} \geq 5^{23/20} 5^{10/100} = 5^{5/4}$$

which contradicts the Fontaine bound. On the other hand, We have the following equality regarding the discriminant ([10], IV. Proposition 4):

$$\sum_{i=0}^{\infty} |G_i| - 1 = v_{\mathfrak{P}}(\mathfrak{D}_{F/E})$$

and so if $v_{F/E}$ is the exponent of the discriminant,

$$v_{F/E} \equiv e_{F/E} - 1 \mod (5-1) \equiv 0 \mod 4.$$

Thus $v_{F/E} = 4$ or $8$. Since we have wild ramification, $v_{F/E} > e_{F/E} - 1$, and thus $v_{F/E} = 8$, and $\Delta_{F/E} = \pi_E^8$.

Suppose now that $E = \mathbb{Q}(\zeta_5, 24^{1/5})$. Let $\pi_{K,i}$ be the unique prime above $\pi_{E,i}$ in $\mathcal{O}_K$. If $\Delta_{L/K} = (\pi_{K,1} \ldots \pi_{K,5})^v$, an argument similar to the above using the Fontaine bound shows that $v < 10$. Thus

$$\Delta_{L/E} = \Delta_{K/E}^5 N_{K/E}(\Delta_{L/K}) < (\pi_{E,1} \ldots \pi_{E,5})^{50+10}.$$

If $\pi_{E,i}$ occurs in $\Delta_{F/E}$ with exponent $v_{F/E}$ then

$$\Delta_{L/E} \geq \Delta_{F/E}^5 = \pi_{E,i}^{5v_{F/E}}$$

and thus $v_{F/E} < 12$. Yet, as above, $v_{F/E} \equiv 0 \mod 4$ and thus $v_{F/E} \leq 8$. □

**Corollary 3.1** *If $F/E$ is ramified, and $F/\mathbb{Q}$ is wildly ramified at $5$ then the conductor $\mathfrak{f}_{E/F}$ is equal to $\pi_E^2$. If $F/\mathbb{Q}$ is tamely ramified, then the conductor divides $(\pi_{E,1} \ldots \pi_{E,5})^2$.*

**Proof**. This follows from the previous lemma, and the conductor-discriminant formula. □

Thus the existence of $F$ will therefore be predicted from the ray class group of $\mathfrak{f}_{F/E}$. We may calculate these groups with the aide of `pari`. The results are tabulated in the table in the appendix (section 5.1), and they indicate the proof is complete, after noting that in all cases when the ray class field is non-trivial, the field $K/E$ is either unramified or has conductor dividing $\mathfrak{f}_{F/E}$.

# 4 $N = 10$.

Let us begin by stating the analogues of theorems in section 2.2.

**Theorem 4.1** *Let $G/\mathbb{Z}[\frac{1}{10}]$ be a finite group scheme of 3-power order such that inertia at $2$ and $5$ acts through a procyclic 3-group. Then $G$ has a filtration by the group schemes $\mathbb{Z}/3\mathbb{Z}$ and $\mu_3$. Moreover, if $G$ is killed by $5$, then $\mathbb{Q}(G) \subseteq H$, where $K := \mathbb{Q}(\sqrt[3]{2}, \sqrt[3]{5}, \zeta_3)$, and $H$ is the Hilbert class field of $K$, which is of degree 3 over $K$.*

**Theorem 4.2 (Schoof)** *Let $p = 2$ or $5$. Let $G/\mathbb{Z}[\frac{1}{p}]$ be a finite group scheme of 3-power order such that inertia at $p$ acts through a procyclic 3-group. Then $G$ has a filtration by the group schemes $\mathbb{Z}/3\mathbb{Z}$ and $\mu_3$. Moreover, the extension group $\mathrm{Ext}^1(\mu_3, \mathbb{Z}/3\mathbb{Z})$ of group schemes over $\mathbb{Z}[\frac{1}{p}]$ is trivial, and there exists an exact sequence of group schemes:*

$$0 \longrightarrow M \longrightarrow G \longrightarrow C \longrightarrow 0$$

*where $M$ is a diagonalizable group scheme over $\mathbb{Z}[\frac{1}{p}]$, and $C$ is a constant group scheme.*

One technical difficulty is that $\mathcal{I}_p(H/\mathbb{Q})$ is not a normal subgroup of $\mathrm{Gal}(H/\mathbb{Q})$, for either $p$ equal 2 or 5. We do however make the following observation: The primes 2 and 5 split into 3 distinct primes in $K$. Moreover, these primes remain inert after passing to $H$. The easiest way to see this is by noting that $H$ is the compositum of $K$ and the Hilbert class field of $\mathbb{Q}(20^{1/3})$. In this field, the primes above 2 and 5 are not principal, and so remain inert in the Hilbert class field. Thus the subgroups $\mathcal{D}_p(H/\mathbb{Q})$ are of index 3 in $\mathrm{Gal}(H/\mathbb{Q})$. Moreover, a natural coset representative for the non-trivial cosets of $\mathrm{Gal}(H/\mathbb{Q})/\mathcal{D}_p(H/\mathbb{Q})$ is given by an element of $\mathcal{I}_{p'}(H/\mathbb{Q})$, where $\{p, p'\} = \{2, 5\}$ as an unordered pair. This leads to the following construction:

**Lemma 4.1** *Let $\{p, p'\} = \{2, 5\}$. Let $\overline{\mathcal{M}} \subset A[3]$ be a $\mathcal{D}_p(H/\mathbb{Q})$ module. Let $\sigma \in \mathcal{I}_{p'}(H/\mathbb{Q})$ be a non-trivial element which does not lie in $\mathcal{D}_p(H/\mathbb{Q})$. Let $\{P_1, P_2, \ldots, P_t\}$ be a generating set for $\overline{\mathcal{M}}$ as a $\mathcal{D}_p(H/\mathbb{Q})$ module. Then*

$$\widehat{\mathcal{M}} = \{P_1, \ldots, P_t, (\sigma - 1)P_1, \ldots, (\sigma - 1)P_t\}$$

*is a $\mathrm{Gal}(\overline{\mathbb{Q}}/\mathbb{Q})$ module.*

By Grothendieck (Theorem 1.2) one finds that as an endomorphism, for $\sigma \in \mathcal{I}_{p'}$, $(\sigma - 1)^2 = 0$ on $A[3]$. Thus $\sigma^2 = 2(\sigma - 1) + 1$, and one sees (since $\mathcal{D}_p(H/\mathbb{Q})$ and $\sigma \in \mathcal{I}_{p'}(H/\mathbb{Q})$ generate $\mathrm{Gal}(H/\mathbb{Q})$) that $\widehat{\mathcal{M}}$ is closed under the action of Galois. □

We now apply this construction not to $\overline{\mathcal{M}}_1(p)$, as in section 2.3, but to $\overline{\mathcal{M}}_2(p)$. Let us assume that $\mathrm{ord}_3(\Phi_{\hat{A}}(p))$ is maximal for some $p \in \{2, 5\}$. If $\kappa = \widehat{\mathcal{M}}_2(p)$, then from Theorem 2.1,

$$\mathrm{ord}_3(\Phi_{\hat{A}'}(p)) - \mathrm{ord}_5(\Phi_{\hat{A}}(p)) = \dim\ \kappa \cap \overline{\mathcal{M}}_2(p) + \dim\ \kappa \cap \overline{\mathcal{M}}_1(p) - \dim\ \kappa.$$

Since $\overline{\mathcal{M}}_2(p) \subseteq \kappa \cap \overline{\mathcal{M}}_1(p)$, we find that this quantity is at least $2t_p - \dim\ \kappa$. On the other hand, from the previous lemma we see that $\dim \kappa \leq 2t_p$, with equality if and only if $\{P_1, \ldots P_t, (\sigma - 1)P_1, \ldots, (\sigma - 1)P_t\}$ are independent inside $A[3]$. Since $\mathrm{ord}_3(\Phi_{\hat{A}}(p))$ is maximal, we have equality. Since the image of $(\sigma-1)$ on $A[3]$ for $\sigma \in \mathcal{I}_{p'}$ is contained within $\overline{\mathcal{M}}_2(p')$ and has dimension at most $t_{p'}$, this immediately proves that $t_p \leq t_{p'}$, and by symmetry, that $t_2 = t_5$. Moreover, equality forces $\overline{\mathcal{M}}_2(p) = \kappa \cap \overline{\mathcal{M}}_1(p)$, and thus by dimension considerations, as a vector space,

$$A[3] = \widehat{\mathcal{M}}_2(p) \oplus \overline{\mathcal{M}}_1(p) \setminus \overline{\mathcal{M}}_2(p).$$

**Lemma 4.2** *For* $\mathrm{ord}_3(\Phi_{\hat{A}}(p))$ *maximal,* $\mathbb{Q}(A[3])$ *is unramified at* $p$.

Consider the decomposition $A[3] = \widehat{\mathcal{M}_2}(p) \oplus \overline{\mathcal{M}_1}(p) \setminus \overline{\mathcal{M}_2}(p)$. By definition, $\mathcal{I}_p$ acts trivially on $\overline{\mathcal{M}_1}(p)$. Thus it suffices to show that $\mathcal{I}_p$ acts trivially on $\widehat{\mathcal{M}_2}(p) = \{P_1, \ldots, (\sigma - 1)P_t\}$. Since $\mathcal{I}_p(H/\mathbb{Q}) = \mathcal{I}_p(H/\mathbb{Q}(\zeta_3))$ for $p \in \{2, 5\}$ we work over this field. Since for $\tau \in \mathcal{I}_p$ the image of $(\tau - 1)$ lies within $\overline{\mathcal{M}_2}(p)$, the action of $\tau \in \mathcal{I}_p$ is represented by a matrix:

$$\tau = \begin{pmatrix} \mathrm{Id}_t & a \\ 0 & \mathrm{Id}_t \end{pmatrix}$$

for some $a \in M_t(\mathbb{F}_p)$. On the other hand, $\widehat{\mathcal{M}_2}(p)$ is a $\mathrm{Gal}(\overline{\mathbb{Q}}/\mathbb{Q})$ module and the action of $\sigma$ is given by:

$$\sigma = \begin{pmatrix} \mathrm{Id}_t & 0 \\ \mathrm{Id}_t & \mathrm{Id}_t \end{pmatrix}.$$

It suffices to prove that $a = 0$, since then we have shown $\mathcal{I}_p$ acts trivially on $A[3]$. Since $\mathbb{Q}(A[3]) \subseteq H$, and since $[H : \mathbb{Q}(\zeta_3)] = 27$, this follows from the following result:

**Sub-lemma 2** *Let* $M \subseteq \mathrm{GL}_2(\mathbb{F}_3[a]/I)$ *be a subgroup of order* 27 *containing the elements:*

$$\left\langle \begin{pmatrix} 1 & a \\ 0 & 1 \end{pmatrix}, \begin{pmatrix} 1 & 0 \\ 1 & 1 \end{pmatrix} \right\rangle$$

*then* $a = 0$.

**Proof.** If $M$ is Abelian, then $[\sigma, \tau] = 1$, and one computes immediately that $a = 0$. From a classification of all non-Abelian groups of order 27, we find that $[M, M]$ is of order 3 and central. $[\sigma, \tau]^3 = 0$ implies that $a^3 = 0$. Assuming this, $[\sigma, \tau]\sigma - \sigma[\sigma, \tau] = 0$ implies that $a^2 = 0$ and $2a + a^2 = 0$. In characteristic 3, this proves that $a = 0$. □

With this result, we may now establish Theorem 2.3 in much the same way as Theorem 2.2. Here are the extra steps required to complete the proof:

1. For $\mathrm{ord}_3(\Phi_{\hat{A}}(2))$ maximal, the exact sequence of group schemes

$$0 \longrightarrow \mu_3^m \longrightarrow A[3] \longrightarrow (\mathbb{Z}/3\mathbb{Z})^n \longrightarrow 0$$

   follows from Theorem 4.2 and Lemma 4.2. Lemma 2.2 applies *mutatis mutandis*.

2. The arguments of section 2.4 and Lemma 2.3 still hold, after noting that $(\mathbb{Z}/10\mathbb{Z})^*$ has order coprime to 3, and thus $\mathbb{Q}$ admits no 3-extension unramified over $\mathbb{Z}[\frac{1}{10}]$.

3. The maximal Galois subextension of $H$ unramified at 2 is $\mathbb{Q}(\zeta_3, 5^{1/3})$. Hence the proof of Lemma 2.4 still applies. Similarly, a proof of Lemma 2.5 requires us only to note that the maximal unramified extension of $\mathbb{Q}(\zeta_3)$ inside $H$ is $\mathbb{Q}(\zeta_3, 10^{1/3})$, which *is* ramified at 5.

4. A final contradiction is reached because

$$3^{4g} \leq (1+\sqrt{3})^{4g}$$

is not true. One might remark at this point that since $A$ has good reduction at 3, and since $A$ is defined over $\mathbb{Q}$, the 3-torsion injects into $A(\mathbb{F}_p)[3]$, as follows from standard facts about formal groups.

Thus is remains to prove Theorem 4.1.

## 4.1 Group Schemes over $\mathbb{Z}[1/10]$.

Since $\text{Gal}(H/\mathbb{Q}(\zeta_3))$ is a 3-group, the discussion at the beginning of section 3 shows that it suffices to prove that if $L = \mathbb{Q}(G \times G_{-1} \times G_2 \times G_5)$ then $L \subseteq H$. One has the following estimate of the root discriminant for $L$:

$$\delta_L < 3^{1+\frac{1}{3-1}} 2^{1-\frac{1}{3}} 5^{1-\frac{1}{3}} = 3^{3/2} 10^{2/3} = 24.118 < 24.258$$

From the estimates of [7] one finds that $[L : \mathbb{Q}] < 280$, and so $[L : K] < 16$. One sees that that $K := \mathbb{Q}(\sqrt{-3}, \sqrt[3]{2}, \sqrt[3]{5}) \subseteq L$. We wish to prove that $\text{Gal}(L/K)$ is a 3-group. The root discriminant of $K$ is $\delta_K = 3^{7/6} 10^{2/3}$, and so $L/K$ is at most ramified at primes above 3. Let $F = \mathbb{Q}(\sqrt{-3}, \sqrt[3]{10})$. Then $F/\mathbb{Q}(\sqrt{-3})$ is unramified at $\sqrt{-3}$. The prime $\sqrt{-3}$ splits completely in $F$, as

$$3 = \pi_1 \pi_2 \pi_3, \qquad \pi_i = \left(3, \left(\frac{\sqrt[3]{10}-1}{\sqrt{-3}}\right) - i\right) \qquad N_{F/\mathbb{Q}}(\pi_i) = 3.$$

The extension $K/F$ is totally ramified at each $\pi_i$, $\pi_i = \mathfrak{p}_i^3$ for all $i$, and $N_{K/\mathbb{Q}}(\mathfrak{p}_i) = 3$.

## 4.2 $L/K$ Tame

In this section we assume that $L/K$ is a *tame* extension.

**Lemma 4.3** $[L : K] \leq 6$.

**Proof**. Arguing as in Lemma 3.2 we find that $N_{K/\mathbb{Q}}(\Delta_{K/\mathbb{Q}}) < 3^{3[L:K]}$. Thus

$$\delta_L = \delta_K N_{K/\mathbb{Q}}(\Delta_{L/K})^{1/[L:\mathbb{Q}]} \leq 3^{7/6} 10^{2/3} 3^{3/18} = 20.082.$$

Yet from the GRH Odlyzko bound, if $[L : \mathbb{Q}] \geq 126$, then $\delta_L > 20.221$. Thus we find $[L : K] \leq 6$. □

If $[L : K] \leq 6$, then either $\text{Gal}(L/K)$ is a 3-group or it surjects onto a non-trivial group of order coprime to 3. In this case, $L$ would contain an Abelian extension $E/K$ tamely ramified and of degree coprime to 3.

**Lemma 4.4** *There are no Abelian extensions $E/K$ tamely ramified of order coprime to 3.*

**Proof.** We proceed via class field theory. According to `pari`, the class number of $K$ is 3, its Hilbert class field being the compositum of $H$ and the Hilbert class field of $\mathbb{Q}(\sqrt{-3}, \sqrt[3]{20})$. Thus it suffices to show that global units of $\mathcal{O}_K$ generate $(\mathcal{O}_K/\mathfrak{p}_1\mathfrak{p}_2\mathfrak{p}_3)^*$. On the other hand, since $K/F$ is totally ramified, we have an isomorphism

$$(\mathcal{O}_K/\mathfrak{p}_1\mathfrak{p}_2\mathfrak{p}_3)^* \simeq (\mathcal{O}_F/\pi_1\pi_2\pi_3)^* \simeq \mathbb{F}_3^* \times \mathbb{F}_3^* \times \mathbb{F}_3^*$$

Hence it suffices to use global units from $\mathcal{O}_F$. Let $v = (\sqrt[3]{10} - 1)/\sqrt{-3}$. Then from `pari` we find that the 2 fundamental units of $\mathcal{O}_F$ are given by

$$\epsilon_1 = \frac{1}{4}v^4 - \frac{1}{2}v^2 + \frac{3}{2}v - \frac{1}{4}$$

$$\epsilon_2 = \frac{1}{4}v^4 - \frac{1}{2}v^3 + \frac{3}{2}v^2 - \frac{1}{4}$$

We find that the images of $-1$, $\epsilon_1$, and $\epsilon_2$ in $\mathcal{O}_F/\pi_1 \times \mathcal{O}_F/\pi_2 \times \mathcal{O}_F/\pi_3$ are $(-1, -1, -1)$, $(1, 1, -1)$ and $(1, -1, 1)$ respectively. Since these elements generate the group $(\mathbb{F}_3^*)^3$, we are done. $\square$

## 4.3 $L/K$ Wild

We assume that $L/K$ is wildly ramified at 3, and (for the moment) not a 3-group. If $\text{Gal}(L/K)^{ab}$ is not a 3-group, then there would exist a corresponding extension $E/K$ tame of order coprime to 3. Since no such extensions exist (see the tame case), we may also assume that $\text{Gal}(L/K)^{ab}$ is a 3-group. Let $G$ denote the group $\text{Gal}(L/K)$. There should be no confusion between the group $G$ and the group scheme $G$, which will not appear again. Let $n = |G|$. Since $n < 16$, $n \in \{6, 12, 15\}$. All groups of order 15 are Abelian. If $n = 6$, the only non-Abelian group is $S_3$. Yet $S_3^{ab} = \mathbb{Z}/2\mathbb{Z}$. Thus $n = 12$. The only group $G$ of order 12 such that $G^{ab} = \mathbb{Z}/3\mathbb{Z}$ is the non-trivial extension of $\mathbb{Z}/3\mathbb{Z}$ by $\mathbb{Z}/2\mathbb{Z} \oplus \mathbb{Z}/2\mathbb{Z}$.

**Lemma 4.5** $N_{L/K}(\Delta_{L/K}) \geq 3^{66}$, $N_{L/K}(\Delta_{L/K}) \leq 3^{69}$.

Assume otherwise. Then since $N_{L/K}(\Delta_{L/K}) = 3^{3k}$ for some $k$, it must be bounded by $3^{63}$. Since $[L:\mathbb{Q}] = 216 = 12 \times 18$,

$$\delta_L = \delta_K N_{K/\mathbb{Q}}(\Delta_{L/K})^{1/[L:\mathbb{Q}]} \leq 3^{7/6} 10^{2/3} 3^{63/216} = 23.039 < 23.089.$$

Yet from the GRH Odlyzko bound, $\delta_L > 23.089$. The other inequality is violated if and only if $N_{L/K}(\Delta_{L/K}) \geq 3^{72}$. Yet in this case,

$$\delta_L = \delta_K N_{K/\mathbb{Q}}(\Delta_{L/K})^{1/[L:\mathbb{Q}]} \geq 3^{7/6} 10^{2/3} 3^{72/216} = 3^{3/2} 10^{2/3}.$$

Yet $\delta_L < 3^{3/2} 10^{2/3}$ by the Fontaine bound, and we are done. $\square$

Before we proceed, we introduce some notation and results from Serre [10]. Let $G_i \subseteq G$ be the decomposition groups of some prime $\mathfrak{p}$ above 3 in $L/K$. These groups are defined by $\mathfrak{p}$ up to conjugacy. However, since $L/\mathbb{Q}$ is Galois, the orders of $G_i$ are independent of the choice of $\mathfrak{p}$ above 3. Let $\mathfrak{P}$ be a prime above $\mathfrak{p}$. Let us simplify some notation. Let $v = v_{\mathfrak{P}}(\mathfrak{D}_{L/K})$, $f = f_{L/K}$, $e = e_{L/K}$, $r = r_{L/K}$, We have equalities:

$$\mathfrak{D}_{L/K} = \prod_{i=1}^{3}\prod_{j=1}^{r}\mathfrak{P}_{i,j}^{v} \qquad \Delta_{L/K} = \prod_{i=1}^{3}\mathfrak{p}_i^{frv} \qquad N_{K/\mathbb{Q}}(\Delta_{L/K}) = 3^{3frv}.$$

From the previous lemma, $22 \leq frv \leq 23$. Moreover, $fre = [L:K] = 12$. Since $K/L$ is wildly ramified, $3|e$. Hence it suffices to show that $e = 3$, $e = 6$ and $e = 12$ all lead to contradictions. If $e = 3$, then $fr = 4$. Yet $fr$ divides 23 or 22, which is impossible. Suppose that $e = 6$. Then $|G_0| = 6$ must be a normal subgroup of $G$ since it is a subgroup of index 2. If $G$ had such a subgroup, then $G^{ab}$ would not be a 3-group. Thus we may assume that $e = 12$. If $e = 12$ then the 3-group $G_1$ would be a normal subgroup of $G_0 = G$. Since $G$ has no such subgroup, we are done, and $\mathrm{Gal}(L/K)$ is a 3-group.

Thus we may assume that $L/K$ is Galois of degree dividing 9, and thus Abelian. Let $\mathfrak{f}_{L/K}$ be the conductor of this extension. If $(\pi_{K,1}\pi_{K,2}\pi_{K,3})^3|\mathfrak{f}_{L/K}$, then from the conductor discriminant formula $\delta_L$ exceeds the Fontaine bound. Thus it suffices to note that that the ray class field of $\mathfrak{f} = (\pi_{K,1}\pi_{K,2}\pi_{K,3})^2$ of $K$ is $\mathbb{Z}/3\mathbb{Z}$, coming exactly from the Hilbert class field $H$ of $K$.

## 4.4 Acknowledgments.

I would like to thank René Schoof for many enlightening discussions.

# 5 Appendix.

## 5.1 Ray Class Fields.

Class field calculations. Here are some computations done using `pari`. They took between 15 minutes and an hour each. The essentials of the `pari` script are below.

| $E$ | $\delta_E$ | $\mathfrak{f}$ | $|\mathrm{Cl}_{\mathfrak{f}}|$ |
|---|---|---|---|
| $\mathbb{Q}(\zeta_5, 2^{1/5})$ | $5^{23/20}2^{4/5}$ | $\pi_E^2$ | 1 |
| $\mathbb{Q}(\zeta_5, 3^{1/5})$ | $5^{23/20}3^{4/5}$ | $\pi_E^2$ | 1 |
| $\mathbb{Q}(\zeta_5, 6^{1/5})$ | $5^{23/20}6^{4/5}$ | $\pi_E^2$ | 5 |
| $\mathbb{Q}(\zeta_5, 12^{1/5})$ | $5^{23/20}6^{4/5}$ | $\pi_E^2$ | 5 |
| $\mathbb{Q}(\zeta_5, 24^{1/5})$ | $5^{3/4}6^{4/5}$ | $(\pi_{E,1}\ldots\pi_{E,5})^2$ | 5 |
| $\mathbb{Q}(\zeta_5, 48^{1/5})$ | $5^{23/20}6^{4/5}$ | $\pi_E^2$ | 5 |
| $\mathbb{Q}(\zeta_3, 2^{1/3}, 5^{1/3})$ | $3^{7/6}10^{2/3}$ | $(\pi_{E,1}\pi_{E,2}\pi_{E,3})^2$ | 3 |

## 5.2 Pari Script.

Here is the `pari` script for fields other than $\mathbb{Q}(\zeta_5, 24^{1/5})$ and $\mathbb{Q}(\zeta_3, 2^{1/3}, 5^{1/3})$, where an adjustment must be made since the conductor is of a slightly different form. The calculation of the discriminant was included as a check against typographical errors in the defining polynomials.

```
allocatemem()
allocatemem()
allocatemem()
allocatemem()
nf=nfinit(poly defining K);
factor(nf[3])
bnf=bnfinit(nf[1],1);
pd=idealprimedec(nf,5);
pd1=idealhnf(nf,pd[1]);
idealnorm(nf,pd1)
pd2=idealmul(nf,pd1,pd1);
bnrclass(bnf,pd1)
bnrclass(bnf,pd2)
```

*Email address*: `fcale@math.berkeley.edu`